\documentclass[12pt,a4paper]{article}
\usepackage{graphicx}
\usepackage[T1]{fontenc}
\usepackage{amsmath,amssymb,anysize,rotating}
\begin{document}
\title{The exponentiated xgammma distribution: Estimation and its application}
\author{Abhimanyu Singh Yadav $^{a}$, Mahendra Saha $^{a}$, Harsh Tripathi$^{a}$\footnote{Corresponding author. e-mail: 2017phdsta01@curaj.ac.in} and Sumit Kumar $^{a}$\\
\small $^{a}$Department of Statistics, Central University of Rajasthan, Rajasthan, India}
\date{}
\maketitle
\begin{abstract}
This article aims to introduced a new lifetime distribution named as exponentiated xgamma distribution (EXGD). The new generalization obtained from xgamma distribution, a special finite mixture of exponential and gamma distributions. The proposed model is very flexible and positively skewed. Different statistical properties of the proposed model, viz., reliability characteristics, moments, generating function, mean deviation, quantile function, conditional moments, order statistics, reliability curves and indices and random variate generation etc. have been derived. The estimation of the of the survival and hazard rate functions of the EXGD has been approached by different methods estimation, viz., moment estimate (ME),maximum likelihood estimate (MLE), ordinary least square and weighted least square estimates (LSE and WLSE), Cram\`er-von-Mises estimate (CME) and maximum product spacing estimate (MPSE). 
At last, one medical data set has been used to illustrate the applicability of the proposed model in real life scenario.  
\end{abstract}
\section{Introduction}
In reliability analysis, the lifetime of any electronic device or items is varying in nature. Hence, it seems to be logical to model the lifetime data with a specific probability distribution. The exponential distribution and it's different generalizations, e.g., Weibull, gamma, exponentiated exponential etc. have been often used to model the data with constant, monotone hazard rate. Also, the finite mixtures of two or more probability distributions are also the better alternative to analyze any life time data, such as Lindley [see, Lineley ($1958$)], generalized Lindley [see, Nadarajah et al.($2011$)]. In the same era of generalization of statistical distributions, the one parameter xgamma distribution (XGD) is one of them, a special finite mixture of exponential and gamma distributions, proposed by Sen et al. ($2016$). In most of the situations, finite mixture distributions arising from the standard distributions play a better role in modelling lifetime phenomena as compared to the standard distributions. Recently, Yadav et al. ($2018$) introduced the inverted version of XGD which possesses the upside-down bathtub-shaped hazard function. The XGD did not provide enough flexibility for analyzing different types of lifetime data as it is of one parameter. It will be useful to consider further alternatives of XGD to increase the flexibility for modelling purposes. In this article, we propose a two parameter family of distribution which generalizes the XGD, named as the exponentiated xgamma distribution (EXGD) and hence the name proposed. The procedure used is based on certain finite mixtures of exponential and gamma distributions. The shape parameter provides more flexibility for describing different types of data allowing hazard rate modelling. Moreover, we also derived some statistical characteristics such moments, conditional moments, order statistics, reliability curves and indices.\\

However, the main objective of this article is two fold: First, we introduced a new probability distribution and studied the several statistical properties of EXGD, as the generalized version of XGD, introduced by Sen et al. ($2016$). Second, different methods of estimation have been employed to estimate the unknown parameters as well as survival function and hazard rate functions of EXGD. To the best of our knowledge thus far, no attempt has been made to introduce the generalized version of XGD. Our aim is to fill up this gap through the present study.\\

Rest of the article is organised as follows: Section $2$ introduced the EXGD and presented its reliability characteristics. Moments, generating function, mean deviation, conditional moments, order statistics and reliability curve etc. haveen discussed in Section $3$. Also we have proposed an algorithm generating random variate from the proposed distribution. Section $4$ discussed the different methods of estimates of the survival function and hazard rate function using maximum likelihood estimate, least square and weighted least square estimate, Cram\`er-von-Mises estimate and maximum product spacing estimate of the parameters. 
 A medical data set is used to illustrate the applicability of the proposed model in real life scenario in Section $5$. Concluding remarks are made in Section $6$.   
\section{The model and its reliability characteristics}
\begin{figure}[htbp]
\centering
\resizebox{12cm}{6cm}{\includegraphics[trim=.01cm 2cm .01cm .01cm]{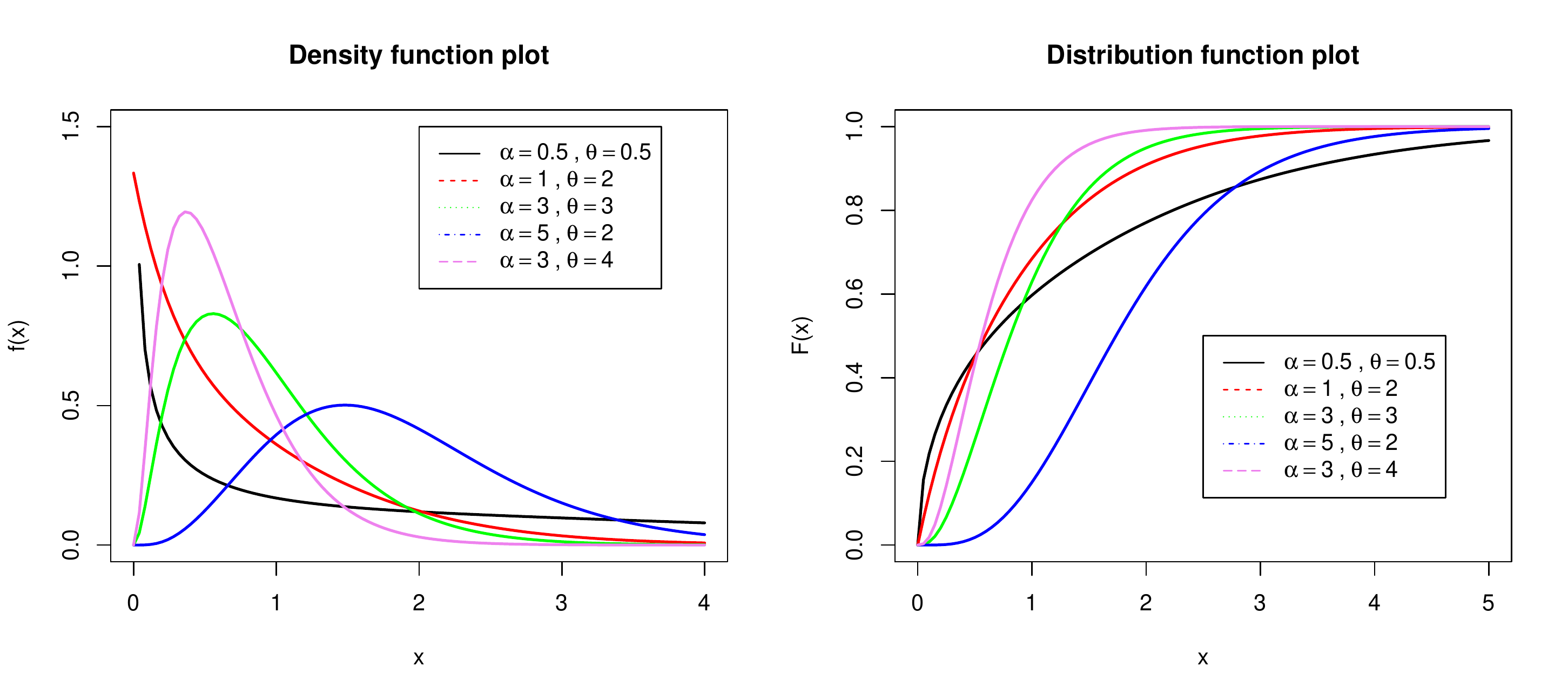}}
\vspace{1cm}
\caption{PDF and CDF plots of EXGD.}
\label{fig1}
\end{figure}

\begin{figure}[htbp]
\centering
\resizebox{12cm}{6cm}{\includegraphics[trim=.01cm 2cm .01cm .01cm]{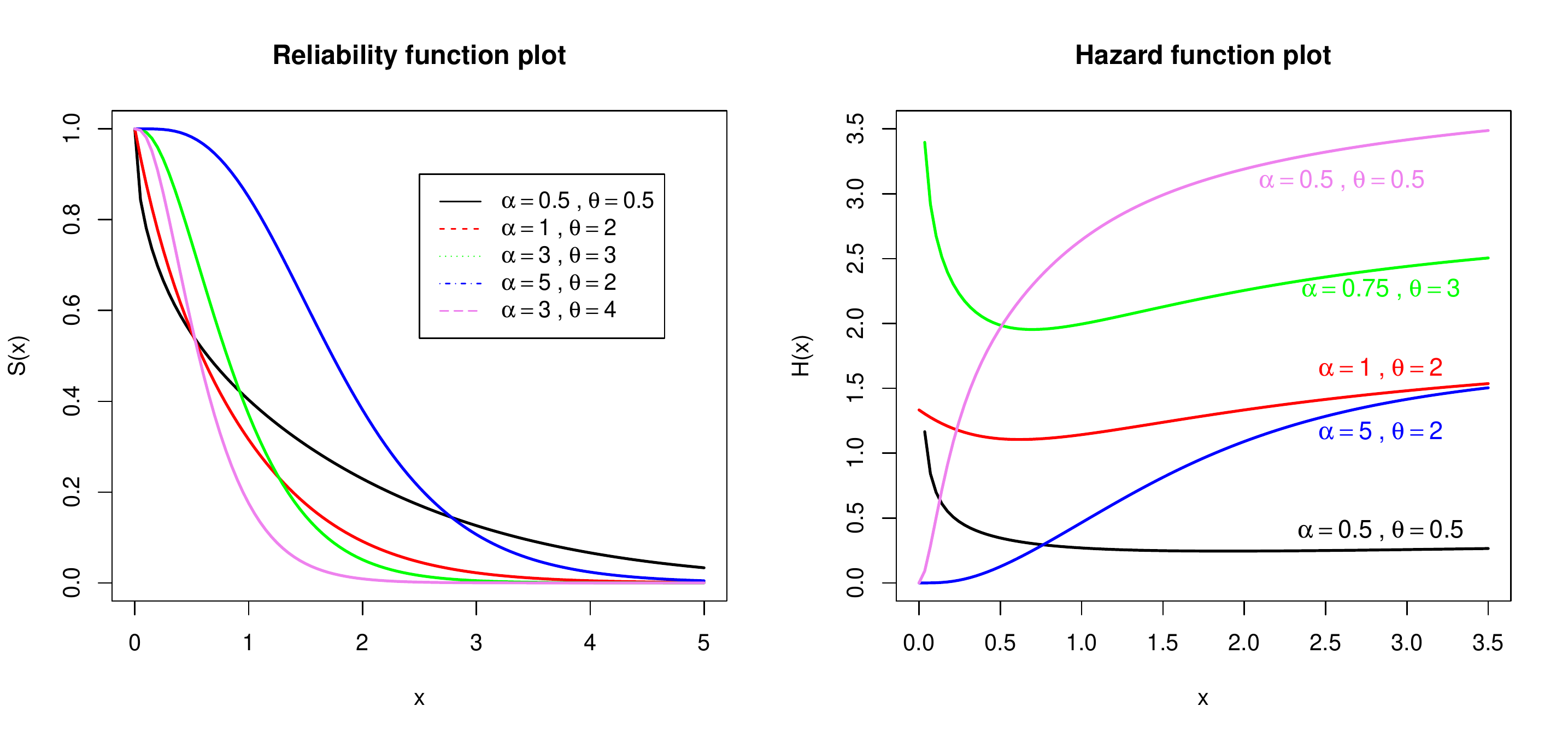}}
\vspace{1cm}
\caption{Survival and hazard rate functions plot EXGD.}
\label{fig2}
\end{figure}

Sen et al. ($2016$) introduced a finite mixture of exponential ($\theta$) and gamma ($3,\theta$) distributions with mixing proportion $\pi_1=\frac{\theta}{1+\theta}$ and $\pi_2=1-\pi_1=\frac{1}{1+\theta}$ to obtain a probability distribution, named as xgamma distribution (XGD), given by the probability density function (PDF) and cumulative distribution function(CDF), respectively as:
\begin{eqnarray}\label{eq1}
f(x;\theta)&=&\frac{\theta^2}{(1+\theta)}\left(1+\frac{\theta}{2}x^2\right)e^{-\theta x}~~;x>0,~\theta>0
\end{eqnarray}
\begin{eqnarray}\label{eq2}
F(x;\theta)&=&1-\frac{\left(1+\theta+\theta x+\frac{\theta^2 x^2}{2}\right)}{(1+\theta)}e^{-\theta x}~~;x>0,~\theta>0
\end{eqnarray}
They have also investigated the mathematical, structural and survival properties and found that in many cases the XGD has more flexibility than the exponential as well as the Lindley distribution. In the era of generalization of new distributions by introducing an extra parameter to any base line distributions, numerous methods are available in literature which accommodated the different shapes of hazard rate. For example, Gupta and Kundu ($1999$) introduced the exponentiated exponential distributions as an alternative to Weibull and gamma distributions, Nadarajah et al. ($2011$) proposed generalized version of the Lindley distribution and shown the superiority of that model rather than the Lindley distribution, exponentiated Rayleigh distribution[see, Surles and Paddgette ($2001$)], exponentiated gamma distribution [see, Shawky and Bakoban ($2006$)], exponentiated Weibull distribution [see, Mudholkar and Srivasthva ($1993$)] etc. All these models are generalized by introducing a shape parameter as power of CDF of the base line model. Obviously, the model with more parameters provides more flexibility but it adds the complexity in the estimation procedure at the same time. Moving on the same path, here we have proposed a new exponentiated distribution, namely, EXGD.\\
Let $X$ be a continuous random variable with CDF ($2$), then the CDF of exponentiated distribution is obtained as

\begin{eqnarray}\label{eq3}
F(x;\alpha,\theta&)=&\left[ 1-\left(1+\theta+\theta x+\frac{\theta^2 x^2}{2}\right)\frac{e^{-\theta x}}{1+\theta}\right]^\alpha;~~~x>0,~\alpha>0,~\theta>0 
\end{eqnarray}
and hence the corresponding PDF is
\begin{eqnarray}\label{eq4}
f(x;\alpha, \theta)&=&\frac{\alpha\theta^2}{1+\theta}  \left[ 1-\left(1+\theta+\theta x+\frac{\theta^2 x^2}{2}\right)\frac{e^{-\theta x}}{1+\theta}\right]^{(\alpha-1)} \left(1+\frac{\theta x^2}{2}\right) e^{-\theta x} 
\end{eqnarray}

where, $\alpha$ is shape parameter and $\theta$ is scale parameter. If shape parameter $\alpha=1$, then Equation($\ref{eq4}$) converted into Equation ($\ref{eq1}$). The shape of PDF and CDF for different values of $\alpha$ and $\theta$ for EXGD are presented in Figure $\ref{fig1}$. PDF plot indicates that EXGD is right-skewed and uni-model or inverted J-shaped.
The basic tools for studying the ageing and associated characteristics of any lifetime equipments are the survival function and hazard rate function. The survival and the hazard rate functions of EXGD($\alpha,\theta$) are given below:
\begin{eqnarray}\label{eq5}
S(x;\alpha,\theta) &=& 1-F(x;\alpha,\theta) \nonumber\\
&=&1-\left[ 1-\left(1+\theta+\theta x+\frac{\theta^2 x^2}{2}\right)\frac{e^{-\theta x}}{1+\theta}\right]^\alpha
\end{eqnarray}
\begin{eqnarray}\label{eq6}
H(x;\alpha,\theta)&=&\left\{\dfrac{\frac{\alpha\theta^2}{1+\theta}  \left[ 1-\left(1+\theta+\theta x+\frac{\theta^2 x^2}{2}\right)\frac{e^{-\theta x}}{1+\theta}\right]^{(\alpha-1)} \left(1+\frac{\theta x^2}{2}\right)e^{-\theta x}}{1-\left[ 1-\left(1+\theta+\theta x+\frac{\theta^2 x^2}{2}\right)\frac{e^{-\theta x}}{1+\theta}\right]^\alpha}\right\} 
\end{eqnarray}
The shapes of the survival and hazard rate functions for EXGD are displayed in Figure $\ref{fig2}$ for certain choices of $\alpha$ and $\theta$. It is to be noted that for $\alpha=1$, the hazard rate function coincide with the hazard rate function of base line distribution [see, Sen et al. ($2016$)], while for other choices of shape and scale parameters EXGD follows the pattern of increasing failure rate (IFR) as well as decreasing failure rate (DFR).
\section{Some statistical properties}
In this section, we have studied some statistical properties EXGD such as moments, generating function, mean deviation, quantile function, conditional moments, order statistics etc. 
\subsection{Moments}
Here, we have provided expression for the moments of the EXGD. The $c$-th order raw moment about origin for EXGD is given as
\begin{equation*}
\begin{split}
E(x^c)&=\int\limits_{0}^{\infty}x^c f(x) dx\\&=\int\limits_{0}^{\infty}x^c \frac{\alpha\theta^2}{1+\theta} e^{-\theta x} \left[ 1-\left(1+\theta+\theta x+\frac{\theta^2 x^2}{2}\right)\frac{e^{-\theta x}}{1+\theta}\right]^{(\alpha-1)} (1+\frac{\theta x^2}{2}) dx\\&=\int\limits_{0}^{\infty}x^c \frac{\alpha\theta^2}{1+\theta} e^{-\theta x} \left[ 1-\left(1+\theta+\theta x+\frac{\theta^2 x^2}{2}\right)\frac{e^{-\theta x}}{1+\theta}\right]^{(\alpha-1)}dx\\&+ \frac{\theta}{2}\int\limits_{0}^{\infty}x^{c+2} \frac{\alpha\theta^2}{1+\theta} e^{-\theta x} \left[ 1-\left(1+\theta+\theta x+\frac{\theta^2 x^2}{2}\right)\frac{e^{-\theta x}}{1+\theta}\right]^{(\alpha-1)} dx
\end{split}
\end{equation*}
Moreover, the expression for moments are not in explicit form, thus, the following results based the lemma $1$, stated below, has been used here to calculate the moments.\\

{\bf Lemma $1$:}\\
Let 
\begin{eqnarray*}
K_1{(a,b,c,\delta)}&=&\int\limits_{0}^{\infty}x^c e^{-\delta x} \left[ 1-\left(1+b+b x+\frac{b^2 x^2}{2}\right)\frac{e^{-b x}}{1+b}\right]^{(a-1)}dx
\\&=&\sum\limits_{i=0}^{\infty}\sum\limits_{j=0}^{i}\sum\limits_{k=0}^{j}\sum\limits_{l=0}^{k}{a-1\choose i}{i\choose j}{j\choose k}{k\choose l}\frac{(-1)^i b^j (b/2)^l \Gamma(c+k+l+1)}{(1+b)^i (bi+\delta)^{(c+k+l+1)}}
\end{eqnarray*}
and 
\begin{eqnarray*}
K_2{(a,b,c,\delta)}&=&\int\limits_{0}^{\infty}x^{c+2} e^{-\delta x} \left[ 1-\left(1+b+b x+\frac{b^2 x^2}{2}\right)\frac{e^{-b x}}{1+b}\right]^{(a-1)}dx
\\&=&\sum\limits_{i=0}^{\infty}\sum\limits_{j=0}^{i}\sum\limits_{k=0}^{j}\sum\limits_{l=0}^{k}{a-1\choose i}{i\choose j}{j\choose k}{k\choose l}\frac{(-1)^i b^j (b/2)^l \Gamma(c+2+k+l+1)}{(1+b)^i (bi+\delta)^{(c+2+k+l+1)}}
\end{eqnarray*}
{\bf Proof:}\\
$$K_1{(a,b,c,\delta)}=\sum\limits_{i=0}^{\infty}{a-1\choose i}\frac{(-1)^i}{(1+b)^i} \int\limits_{0}^{\infty}x^c e^{-\delta x-ibx}\left(1+b+bx+\frac{b^2 x^2}{2}\right)^i dx$$
$$=\sum\limits_{i=0}^{\infty}{a-1\choose i}\frac{(-1)^i}{(1+b)^i} \sum\limits_{j=0}^{i}{i\choose j}b^j \sum\limits_{k=0}^{j}{j\choose k}\int\limits_{0}^{\infty} x^{(c+k)} {\left(1+\frac{bx}{2}\right)}^k e^{-\delta x-ibx} dx$$
\begin{eqnarray}
=\sum\limits_{i=0}^{\infty}{a-1\choose i}\frac{(-1)^i}{(1+b)^i} \sum\limits_{j=0}^{i}{i\choose j}b^j \sum\limits_{k=0}^{j}{j\choose k} \sum\limits_{l=0}^{k}{k\choose l}(b/2)^l \int\limits_{0}^{\infty} x^{c+k+l}e^{-\delta x-ibx} dx
\end{eqnarray}
by use of gamma function, the above equation is written as
$$
K_1(a,b,c,\delta)=\sum\limits_{i=0}^{\infty}\sum\limits_{j=0}^{i}\sum\limits_{k=0}^{j}\sum\limits_{l=0}^{k}{a-1\choose i}{i\choose j}{j\choose k}{k\choose l}\frac{(-1)^i b^j (b/2)^l \Gamma(c+k+l+1)}{(1+b)^i (bi+\delta)^{(c+k+l+1)}}
$$
Similarly solve the $K_2(a,b,c,\delta)$-
$$
K_2(a,b,c,\delta)=\sum\limits_{i=0}^{\infty}\sum\limits_{j=0}^{i}\sum\limits_{k=0}^{j}\sum\limits_{l=0}^{k}{a-1\choose i}{i\choose j}{j\choose k}{k\choose l}\frac{(-1)^i b^j (b/2)^l \Gamma(c+2+k+l+1)}{(1+b)^i (bi+\delta)^{(c+2+k+l+1)}}
$$
By putting $\alpha=a$, $\theta=b$, $\delta=\theta$, $c=r$ in the above lemma $1$, the expression of $r$-th raw moment is given as
\begin{eqnarray}
E(x^r)=\frac{\alpha\theta^2}{1+\theta} \left[K_1(\alpha,\theta,r,\theta)+\frac{\theta}{2} K_2(\alpha,\theta,r,\theta)\right]
\end{eqnarray}
Hence, the first four raw moments of EXGD are obtained as
$$E(x)=\frac{\alpha\theta^2}{1+\theta} \left[K_1(\alpha,\theta,1,\theta)+\frac{\theta}{2} K_2(\alpha,\theta,1,\theta)\right]$$
$$E(x^2)=\frac{\alpha\theta^2}{1+\theta} \left[K_1(\alpha,\theta,2,\theta)+\frac{\theta}{2} K_2(\alpha,\theta,2,\theta)\right]$$
$$E(x^3)=\frac{\alpha\theta^2}{1+\theta} \left[K_1(\alpha,\theta,3,\theta)+\frac{\theta}{2} K_2(\alpha,\theta,3,\theta)\right]$$
$$E(x^4)=\frac{\alpha\theta^2}{1+\theta} \left[K_1(\alpha,\theta,4,\theta)+\frac{\theta}{2} K_2(\alpha,\theta,4,\theta)\right]$$
Again, the first four central moments can be obtained by using the relation between the raw moments and central moments. Hence, the Pearson measures of skewness (SK) and kurtosis (KR) based on second($\mu_2$), third ($\mu_3$) and fourth ($\mu_4$) central moments are obtained by using the following relation, given below:
$$
SK=\frac{\mu_3^2}{\mu_2^3}~~\mbox{and}~~KR= \frac{\mu_4}{\mu_2^2}
$$
\subsection{Generating functions}
Here, in this subsection, the different generating functions, namely, moment generating function $M_x(t)$, characteristics function $\Phi_x(t)$ and Kumulants generating function $K_x(t)$ are derived and presented in the following equations: 
\begin{eqnarray}
M_x(t)=E(e^{tx})&=&\nonumber \int\limits_{0}^{\infty} e^{tx} f(x)dx\\&=&\int\limits_{0}^{\infty} \frac{\theta^2}{(1+\theta)}\left(1+\frac{\theta}{2}.x^2\right)e^{-x(\theta-t)} dx
\end{eqnarray}
By using the same lemma, the moment generating is given as
$$
M_x(t)=\frac{\alpha\theta^2}{1+\theta} \left[K_1(\alpha,\theta,0,\theta-t)+\frac{\theta}{2} K_2(\alpha,\theta,0,\theta-t)\right]
$$
The characteristic function for EXGD is simply obtained by replacing dummy parameter $t$ by $it$, where, $i^2=-1$, given as
$$\Phi_x(t)=\frac{\alpha\theta^2}{1+\theta} \left[K_1(\alpha,\theta,0,\theta-it)+\frac{\theta}{2} K_2(\alpha,\theta,0,\theta-it)\right]$$
The kumulants generating function is the logarithm of the moment generating function and is obtained as
$$
K_x(t)=\log \Phi_x(t)=\log\left(\frac{\alpha\theta^2}{1+\theta}\right)+\log\left[K_1(\alpha,\theta,0,\theta-it)+\frac{\theta}{2} K_2(\alpha,\theta,0,\theta-it)\right]
$$
\subsection{Mean deviation}
The mean deviation about mean of random variable $X$, having density function ($3$) is obtained by;
$$
M.D=\int\limits_{0}^{\infty}|(x-\mu)| f(x)  dx
$$
where, $\mu$=E(x) 
$$
M.D=\int\limits_{o}^{\mu}(\mu-x)f(x)dx+\int\limits_{\mu}^{\infty}(x-\mu)f(x)dx
$$
$$
M.D=\mu F(\mu)-\int\limits_{0}^{\mu}x f(x)dx+\int\limits_{\mu}^{\infty}x f(x) dx-\mu+\mu F(\mu)
$$
$$
M.D=2\mu F(\mu)-2\mu+2\int\limits_{\mu}^{\infty}x f(x)dx
$$
where, $F(\mu)$ stands for CDF of X upto point $\mu$ and \\
$$
\int\limits_{\mu}^{\infty}x f(x) dx=\frac{\alpha\theta^2}{1+\theta}\left[L_1(\alpha,\theta,1,\theta,x)+\frac{\theta}{2}L_2(\alpha,\theta,1,\theta,x)\right]
$$
Using the value from the above integral one can evaluate mean deviation about mean. 
\subsection{Quantile Function }
If $Q(p)$ be the quantile of order $p$ of the EXGD random variable $X$, then the quantile function will be the solution of the following equation
\begin{eqnarray}
p&=&\left[ 1-\left(1+\theta+\theta Q(p)+\frac{\theta^2 Q(p)^2}{2}\right)\frac{e^{-\theta Q(p)}}{1+\theta}\right]^\alpha
\end{eqnarray}
The skewness and kurtosis are the two important measures to study the symmetry and convexity of the curve. The Bowley measure of skewness [see, Bowley ($1920$)] and Moors measure of kurtosis [see, Moors ($1988$)] based on quantile can be used and are given as follows:
\begin{eqnarray*}
SK&=&\frac{Q(\frac{3}{4})-2Q(\frac{1}{2})+Q(\frac{1}{4})}{Q(\frac{3}{4})-Q(\frac{1}{4})}
\end{eqnarray*}
\begin{eqnarray*}
KR&=&\frac{Q(\frac{7}{8})-Q(\frac{5}{8})+Q(\frac{3}{8})-Q(\frac{1}{8})}{Q(\frac{6}{8})-Q(\frac{2}{8})}
\end{eqnarray*}

\subsection{Conditional moments}
The conditional moments about origin is defined as;
$$E(X^n|X>x)=\int\limits_{x}^{\infty} x^n \frac{f(x)}{1-F(x)}$$
where, F(x) is CDF of EXGD, given in Equation ($3$).
\begin{equation}
\begin{split}
E(X^n|X>x)&=\frac{1}{1-F(x)} \frac{\alpha\theta^2}{1+\theta}\left[\int\limits_{x}^{\infty} x^n \eta(x_i,\theta)dx+ \int\limits_{x}^{\infty} \frac{\theta}{2}x^{n+2}\eta(x_i,\theta) dx\right]
\end{split}
\end{equation}
where; $\eta(x_i,\theta)= e^{-\theta x} \left[ 1-\left(1+\theta+\theta x+\frac{\theta^2 x^2}{2}\right)\frac{e^{-\theta x}}{1+\theta}\right]^{(\alpha-1)}$. 
The above Equation ($11$) involved two integral which can not be easily tractable. Thus, the following lemma is used to evaluate the integral.\\ 

{\bf Lemma $2$:}\\
Let
\begin{eqnarray*}
L_1(a,b,c,\delta,t)&=&\int\limits_{t}^{\infty} x^c e^{-\delta x} \left[ 1-\left(1+b+b x+\frac{b^2 x^2}{2}\right)\frac{e^{-b x}}{1+b}\right]^{(a-1)} dx
\\&=&\sum\limits_{i=0}^{\infty}\sum\limits_{j=0}^{i}\sum\limits_{k=0}^{j}\sum\limits_{l=0}^{k}{a-1\choose i}{i\choose j}{j\choose k}{k\choose l}\frac{(-1)^i b^j (b/2)^l}{(1+b)^i} \frac{\Gamma{c+k+l+1,t(bi+\delta)}}{(bi+\delta)^{(c+k+l+1)}}
\end{eqnarray*}
and
\begin{eqnarray*}
L_2(a,b,c,\delta,t&)=&\int\limits_{t}^{\infty} x^{c+2} e^{-\delta x} \left[ 1-\left(1+b+b x+\frac{b^2 x^2}{2}\right)\frac{e^{-b x}}{1+b}\right]^{(a-1)} dx
\\&=&\sum\limits_{i=0}^{\infty}\sum\limits_{j=0}^{i}\sum\limits_{k=0}^{j}\sum\limits_{l=0}^{k}{a-1\choose i}{i\choose j}{j\choose k}{k\choose l}\frac{(-1)^i b^j (b/2)^l}{(1+b)^i} \frac{\Gamma{c+2+k+l+1,t(bi+\delta)}}{(bi+\delta)^{(c+2+k+l+1)}}
\end{eqnarray*}
{\bf Proof:} The proof of the above lemmas are straight forward as the previous one.\\
$$
L_1(a,b,c,\delta,t)=\sum\limits_{i=0}^{\infty}\sum\limits_{j=0}^{i}\sum\limits_{k=0}^{j}\sum\limits_{l=0}^{k}{a-1\choose i}{i\choose j}{j\choose k}{k\choose l}\frac{(-1)^i b^j (b/2)^l}{(1+b)^i} \int\limits_{t}^{\infty} x^{c+k+l} e^{-x(bi+\delta)} dx 
$$
After simplifications,
$$
L_1(a,b,c,\delta,t)=\sum\limits_{i=0}^{\infty}\sum\limits_{j=0}^{i}\sum\limits_{k=0}^{j}\sum\limits_{l=0}^{k}{a-1\choose i}{i\choose j}{j\choose k}{k\choose l}\frac{(-1)^i b^j (b/2)^l}{(1+b)^i} \frac{\Gamma{c+k+l+1,t(bi+\delta)}}{(bi+\delta)^{(c+k+l+1)}} 
$$
similarly 
$$
L_2(a,b,c,\delta,t)=\sum\limits_{i=0}^{\infty}\sum\limits_{j=0}^{i}\sum\limits_{k=0}^{j}\sum\limits_{l=0}^{k}{a-1\choose i}{i\choose j}{j\choose k}{k\choose l}\frac{(-1)^i b^j (b/2)^l}{(1+b)^i} \frac{\Gamma{c+2+k+l+1,t(bi+\delta)}}{(bi+\delta)^{(c+2+k+l+1)}} 
$$
Expression of $E(X^n|X>x)$ given below:
$$E(X^n|X>x)=\frac{1}{1-F(x)} \frac{\alpha\theta^2}{1+\theta}\left[L_1(\alpha,\theta,n,\theta,x)+\frac{\theta}{2}L_2(\alpha,\theta,n,\theta,x)\right]$$
Using Lemma $2$, the conditional moments are given as
$$E(X|X>x)=\frac{1}{1-F(x)} \frac{\alpha\theta^2}{1+\theta}\left[L_1(\alpha,\theta,1,\theta,x)+\frac{\theta}{2}L_2(\alpha,\theta,1,\theta,x)\right]$$
$$E(X^2|X>x)=\frac{1}{1-F(x)} \frac{\alpha\theta^2}{1+\theta}\left[L_1(\alpha,\theta,2,\theta,x)+\frac{\theta}{2}L_2(\alpha,\theta,2,\theta,x)\right]$$
$$E(X^3|X>x)=\frac{1}{1-F(x)} \frac{\alpha\theta^2}{1+\theta}\left[L_1(\alpha,\theta,3,\theta,x)+\frac{\theta}{2}L_2(\alpha,\theta,3,\theta,x)\right]$$
$$E(X^4|xX>x)=\frac{1}{1-F(x)} \frac{\alpha\theta^2}{1+\theta}\left[L_1(\alpha,\theta,4,\theta,x)+\frac{\theta}{2}L_2(\alpha,\theta,4,\theta,x)\right]$$

\subsection{Order statistics}
Let $X_1,~X_2,~X_3,...,~X_n$ is a random sample of size from EXGD. Then, the ordered observations as $X_{(1)}< X_{(2)}<X_{(3)}<.....<X_{(n)}$ constitute the order statistic. Let $X_{(k;n)}$ denotes the $k$-th order statistic, then the PDF and CDF of $k$-th order statistic are computed as
$$
f(X_{(k;n)}=t)=\frac{n !}{(n-k)! (k-1)!}\sum\limits_{l=0}^{n-k}{n-k\choose l}{(-1)^l}~{[F(t)]^l}~F^{k-1}(t)~f(t)
$$
Using Equation ($3$) and Equation ($4$), the PDF of $X_{(k;n)}$ is 
{\scriptsize 
\begin{eqnarray}
f(X_{(k;n)}=t)=\frac{\alpha\theta^2}{1+\theta} \frac{n !}{(n-k)! (k-1)!} e^{-\theta t} (1+\frac{\theta t^2}{2}) \sum\limits_{l=0}^{n-k}{n-k\choose l}{(-1)^l}{\left[ 1-\left(1+\theta+\theta t+\frac{\theta^2 t^2}{2}\right)\frac{e^{-\theta t}}{1+\theta}\right]}^{\alpha(k+l)-1}
\end{eqnarray}}
The CDF of $k$-th order statistic is
$$
F(X_{(k;n)}=t)=\sum\limits_{j=k}^{n} \sum\limits_{l=0}^{n-j} {n\choose j} {n-j\choose l} {(-1)^l} {F^{(j+l)}(t)}
$$
Using Equation ($3$)
\begin{eqnarray}
F(X_{(k;n)}=t)=\sum\limits_{j=k}^{n} \sum\limits_{l=0}^{n-j} {n\choose j} {n-j\choose l} {(-1)^l} {\left[ 1-\left(1+\theta+\theta t+\frac{\theta^2 t^2}{2}\right)\frac{e^{-\theta t}}{1+\theta}\right]}^{\alpha(j+l)}
\end{eqnarray}
The distribution $X_{(1)}=\min(X_{(1)}< X_{(2)}<X_{(3)}<.....<X_{(n)})$ and $X_{(n)}=\max(X_{(1)}< X_{(2)}<X_{(3)}<.....<X_{(n)})$ can be computed with help of above Equations by putting $k=1$ and $k=n$ respectively. 
\subsection{Reliability curves and indices}

Bonferroni and Lorenz curves are very important tools in actuarial and population science to study the income and poverty level. Besides these filed, the reliability curve also evaluated based on specific probability distributions. Let $X$ be a random variable with PDF $f(x)$, defined in Equation ($4$) then Bonferroni curve $B(p)$ and Lorenz curve $L(p)$ are defined by the following Equations ($14$) and ($15$)  
$$B(p)=\frac{1}{p\mu} \int\limits_{0}^{q} x f(x) dx$$
$$B(p)=\frac{1}{p\mu}{[\mu-\int\limits_{q}^{\infty}x f(x) dx]}$$
\begin{eqnarray}
B(p)=\frac{1}{p\mu} \left[\mu-\frac{1}{1-F(x)} \frac{\alpha\theta^2}{1+\theta}\left(L_1(\alpha,\theta,1,\theta,x)+\frac{\theta}{2}L_2(\alpha,\theta,1,\theta,x)\right)\right]
\end{eqnarray}
and
$$L(p)=\frac{1}{\mu} \int\limits_{0}^{q} x f(x) dx$$
$$L(p)=\frac{1}{\mu}{[\mu-\int\limits_{q}^{\infty} x f(x) dx]}$$
\begin{eqnarray}
L(p)=\frac{1}{\mu} \left[\mu-\frac{1}{1-F(x)} \frac{\alpha\theta^2}{1+\theta} \left(L_1(\alpha,\theta,1,\theta,x)+\frac{\theta}{2}L_2(\alpha,\theta,1,\theta,x)\right)\right]
\end{eqnarray}
where, $\mu=E(x)$ and the indices based on these two curves are given as
$$B=1-\int\limits_{0}^{1} B(p) dp$$
$$G=1-2\int\limits_{0}^{1} L(p)dp$$

\subsection{Random number generation}
To generate random number from EXGD with shape parameter $(\alpha)$ and scale parameter $(\theta)$. The following steps may be used.\\
\begin{enumerate}
\item Specified the values of $\alpha$, $\theta$ and $n$.
\item Generate $U_i$ from uniform(0,1) distribution $(i=1,2,3...,n)$.\\
\item Generate $V_i$ from $gamma(\alpha,\theta)$ distribution $(i=1,2,3...,n)$.\\
\item Generate $W_i$ from $gamma(\alpha+2,\theta)$ distribution $(i=1,2,3...,n)$.\\
\item  If $U_i\leq\frac{\theta}{\theta+1}$, set $X_i=V_i$, otherwise set $X_i=W_i$.\\
\end{enumerate}
If we take $\alpha=1$, then we get the random variates from XGD.
\section{Different methods of estimation} 
In this section, we have used five methods of estimation to estimate the unknown parameters as well as survival function $S(t)$ and hazard rate function $H(t)$, namely, moment estimation (ME), maximum likelihood estimation (MLE), least squares estimation (LSE), weighted least squares estimation (WLSE), Cram\`er-von-Mises estimator estimation (CME) and maximum product of spacings estimation (MPSE) respectively for the EXGD.\\

\subsection{Maximum likelihood estimator}
Let $X_{1},~X_{2}, \ldots,~X_{n}$ be a random sample of size $n$ from Equation (\ref{eq4}). Then, the log-likelihood function for the observed random sample $x_{1},~x_{2}, \ldots,~x_{n}$ is given as
\begin{equation}\label{eq19}
\ell(\alpha, \theta)=n\log\alpha+2n\log\theta-n\log(1+\theta)-\theta\sum\limits_{i=1}^{n}x_i+(\alpha-1)\sum\limits_{i=1}^{n}\log U(x_i)+\sum\limits_{i=1}^{n}\log\left(1+\frac{\theta x_i^2}{2}\right),
\end{equation}
where,
$$
U(x_i)=\left[ 1-\left(1+\theta+\theta x_i+\frac{\theta^2 x_i^2}{2}\right)\frac{e^{-\theta x_i}}{1+\theta}\right].
$$
The resulting partial derivatives of the log-likelihood function are
\begin{equation}\label{eq20}
\frac{\partial\,\ell(\alpha, \theta)}{\partial\alpha}=\frac{n}{\alpha}+\sum\limits_{i=1}^{n}\log U(x_i)=0
\end{equation}
\begin{eqnarray}\label{eq21}
\frac{\partial\,\ell(\alpha, \theta)}{\partial\theta}&=&-\frac{n(\theta+2)}{\theta(1+\theta)}-\sum\limits_{i=1}^{n}\frac{(x_i^2/2)}{(1+{(\theta x_i^2/2)})}-\sum\limits_{i=1}^{n} x_i\nonumber \\&&+\frac{(1-\alpha)}{(1+\theta)^2}\sum\limits_{i=1}^{n} \frac{e^{-\theta x_i}(2\theta x_i+\theta^2 x_i+{(\theta^2 x_i^2/2)}+{(\theta^2 x_i^3/2)}+{(\theta^3 x_i^3/2)})}{U(X_i)}
\end{eqnarray}
Equating these partial derivatives to zero does not yield closed-form solutions for the MLEs and thus a numerical method is used for solving these equations simultaneously. Substituting the MLEs, we can get the estimators of $S(x)$ and $H(x)$, given as
\begin{eqnarray}\label{eq22}
\hat{S}(x)_{mle}&=& 1-\left[ 1-\left(1+{\hat{\theta}_{mle}}+{\hat{\theta}_{mle}} x+\frac{{\hat{\theta}_{mle}}^2 x^2}{2}\right)\frac{e^{-{\hat{\theta}_{mle}} x}}{1+{\hat{\theta}_{mle}}}\right]^{\hat{\alpha}_{mle}}
\end{eqnarray}
and
\begin{eqnarray}\label{eq23}
\hat{H}(x)_{mle}&=&\left\{\dfrac{\frac{\alpha\hat{\theta}_{mle}^2}{1+\hat{\theta}_{mle}}  \left[ 1-\left(1+\hat{\theta}_{mle}+\hat{\theta}_{mle} x+\frac{\hat{\theta}_{mle}^2 x^2}{2}\right)\frac{e^{-\hat{\theta}_{mle} x}}{1+\hat{\theta}_{mle}}\right]^{({\hat{\alpha}_{mle}}-1)} \left(1+\frac{\hat{\theta}_{mle} x^2}{2}\right)e^{-\hat{\theta}_{mle} x}}{1-\left[ 1-\left(1+\hat{\theta}_{mle}+\hat{\theta}_{mle} x+\frac{\hat{\theta}_{mle}^2 x^2}{2}\right)\frac{e^{-\hat{\theta}_{mle} x}}{1+\hat{\theta}_{mle}}\right]^{\hat{\alpha}_{mle}}}\right\} 
\end{eqnarray}
\subsection{Least square and weighted least square estimator}
The least square estimator (LSE) and the weighted least square estimator (WLSE) were proposed by Swain et al. ($1988$) to estimate the parameters of Beta distribution. Suppose $F(x_{(j)})$ denotes the distribution function of the ordered random variables $x_{(1)}<x_{(2)}<\cdots <x_{(n)}$, where, $\{x_{1},x_{2},\cdots ,x_{n}\}$ is a random sample of size $n$ from a distribution function $F(\cdot)$. Therefore, in this case, the LSE of $\alpha$, $\theta$, say, $\hat{\alpha}_{lse}$, $\hat{\theta}_{lse}$ can be obtained by minimizing

$$
\mathcal L(\alpha,\theta)=\sum\limits_{i=1}^{n} \left[F(x_{i:n}|\alpha,\theta)-\frac{i}{n+1}\right]^2
$$
 with respect to $\alpha$ and $\theta$, where, $F(\cdot)$ is the CDF, given in Equation (\ref{eq3}). Equivalently, it can be obtained by solving:
$$
\frac{\partial \mathcal L}{\partial \alpha}= \sum\limits_{i=1}^{n} \left[F(x_{i:n}|\alpha,\theta)-\frac{i}{n+1}\right] \frac{\partial F(x_{i:n}|\alpha,\theta)}{\partial \alpha}=0
$$ 
$$
\frac{\partial \mathcal L}{\partial \theta}= \sum\limits_{i=1}^{n} \left[F(x_{i:n}|\alpha,\theta)-\frac{i}{n+1}\right] \frac{\partial F(x_{i:n}|\alpha,\theta)}{\partial \theta}=0
$$
where, 
\begin{eqnarray}
\frac{\partial F(x;\alpha,\theta)}{\partial \alpha}=\left[ 1-\left(1+\theta+\theta x+\frac{\theta^2 x^2}{2}\right)\frac{e^{-\theta x}}{1+\theta}\right]^\alpha \log\left[ 1-\left(1+\theta+\theta x+\frac{\theta^2 x^2}{2}\right)\frac{e^{-\theta x}}{1+\theta}\right]
\end{eqnarray}
and
\begin{eqnarray}
\frac{\partial F(x;\alpha,\theta)}{\partial \theta}&=&\alpha \left[ 1-\left(1+\theta+\theta x+\frac{\theta^2 x^2}{2}\right)\frac{e^{-\theta x}}{1+\theta}\right]^{\alpha-1}\nonumber\\&& \left[\frac{e^{-\theta x}}{(1+\theta)^2}{(1+\theta+\theta x+\frac{\theta^2 x^2}{2})}{(1+x+\theta x)}-{(1+x+\theta x^2)\frac{e^{-\theta x}}{(1+\theta)}}\right]
\end{eqnarray}
Hence, the estimated survival and hazard rate functions are obtained as
\begin{eqnarray}\label{eq22}
\hat{S}(x)_{lse}&=& 1-\left[ 1-\left(1+{\hat{\theta}_{lse}}+{\hat{\theta}_{lse}} x+\frac{{\hat{\theta}_{lse}}^2 x^2}{2}\right)\frac{e^{-{\hat{\theta}_{lse}} x}}{1+{\hat{\theta}_{lse}}}\right]^{\hat{\alpha}_{lse}}
\end{eqnarray}
and
\begin{eqnarray}\label{eq23}
\hat{H}(x)_{lse}&=&\left\{\dfrac{\frac{\alpha\hat{\theta}_{lse}^2}{1+\hat{\theta}_{lse}}  \left[ 1-\left(1+\hat{\theta}_{lse}+\hat{\theta}_{lse} x+\frac{\hat{\theta}_{lse}^2 x^2}{2}\right)\frac{e^{-\hat{\theta}_{lse} x}}{1+\hat{\theta}_{lse}}\right]^{({\hat{\alpha}_{lse}}-1)} \left(1+\frac{\hat{\theta}_{lse} x^2}{2}\right)e^{-\hat{\theta}_{lse} x}}{1-\left[ 1-\left(1+\hat{\theta}_{lse}+\hat{\theta}_{lse} x+\frac{\hat{\theta}_{lse}^2 x^2}{2}\right)\frac{e^{-\hat{\theta}_{lse} x}}{1+\hat{\theta}_{lse}}\right]^{\hat{\alpha}_{lse}}}\right\} 
\end{eqnarray}

WLSE proposed by Swain et al.(1988). Suppose $x_{(1)}, x_{(1)},...x_{(n)}$ be the ordered sample from EXGD. Then WLSE of $\alpha$ and $\theta$ to be estimated by minimizing the function $W(\alpha,\theta)$\\
$$
\mathcal W(\alpha,\theta)=\sum\limits_{i=1}^{n} \frac{{(n+1)^2}{(n+2)}}{i{(n-i+1)}}\left[F(x_{i:n}|\alpha,\theta)-\frac{i}{n+1}\right]^2
$$
Partially derivative of $W(\alpha,\theta)$ with respect to $\alpha$ and $\theta$ and equating both equation equal to zero\\
$$
\frac{\partial \mathcal W}{\partial \alpha}= \sum\limits_{i=1}^{n} \frac{{(n+1)^2}{(n+2)}}{i{(n-i+1)}} \left[F(x_{i:n}|\alpha,\theta)-\frac{i}{n+1}\right] \frac{\partial F(x_{i:n}|\alpha,\theta)}{\partial \alpha}=0$$ 
$$
\frac{\partial \mathcal W}{\partial \theta}= \sum\limits_{i=1}^{n} \frac{{(n+1)^2}{(n+2)}}{i{(n-i+1)}} \left[F(x_{i:n}|\alpha,\theta)-\frac{i}{n+1}\right] \frac{\partial F(x_{i:n}|\alpha,\theta)}{\partial \theta}=0
$$ 
where, $\frac{\partial F(x_{i:n}|\alpha,\theta)}{\partial \alpha}$, $\frac{\partial F(x_{i:n}|\alpha,\theta)}{\partial \theta}$ is defined in Equations ($21$), ($22$) respectively. Hence, the estimated survival and hazard rate functions are obtained as
\begin{eqnarray}\label{eq22}
\hat{S}(x)_{wlse}&=& 1-\left[ 1-\left(1+{\hat{\theta}_{wlse}}+{\hat{\theta}_{wlse}} x+\frac{{\hat{\theta}_{wlse}}^2 x^2}{2}\right)\frac{e^{-{\hat{\theta}_{wlse}} x}}{1+{\hat{\theta}_{wlse}}}\right]^{\hat{\alpha}_{wlse}}
\end{eqnarray}
and
\begin{eqnarray}\label{eq23}
\hat{H}(x)_{wlse}&=&\left\{\dfrac{\frac{\alpha\hat{\theta}_{wlse}^2}{1+\hat{\theta}_{wlse}}  \left[ 1-\left(1+\hat{\theta}_{wlse}+\hat{\theta}_{wlse} x+\frac{\hat{\theta}_{wlse}^2 x^2}{2}\right)\frac{e^{-\hat{\theta}_{wlse} x}}{1+\hat{\theta}_{wlse}}\right]^{({\hat{\alpha}_{wlse}}-1)} \left(1+\frac{\hat{\theta}_{wlse} x^2}{2}\right)e^{-\hat{\theta}_{wlse} x}}{1-\left[ 1-\left(1+\hat{\theta}_{wlse}+\hat{\theta}_{wlse} x+\frac{\hat{\theta}_{wlse}^2 x^2}{2}\right)\frac{e^{-\hat{\theta}_{wlse} x}}{1+\hat{\theta}_{wlse}}\right]^{\hat{\alpha}_{wlse}}}\right\} 
\end{eqnarray}

\subsection{Cramer-von-mises estimator}
To motivate our choice of Cramer-von Mises type minimum distance estimators, MacDonald ($1971$) provided empirical evidence that the bias of the estimator is smaller than the other minimum distance estimators. Thus, The Cramer-von Mises estimators of $\alpha$ and $\theta$, say, $\hat{\alpha}_{cme}$ $\hat{\theta}_{cme}$ receptively can be obtained by minimizing the function 
$$
\mathcal C(\alpha,\theta)=\frac{1}{12n}+\sum\limits_{i=1}^{n}{\left(F(x_{i:n}|\alpha,\theta)-\frac{2i-1}{2n}\right)^2}
$$
with respect to $\alpha$ and $\theta$. The estimator can also be obtained by solving the non linear equations

$$
\frac{\partial \mathcal C}{\partial \alpha}=\sum\limits_{i=1}^{n}{\left(F(x_{i:n}|\alpha,\theta)-\frac{2i-1}{2n}\right)}\frac{\partial F(x_{i:n}|\alpha,\theta)}{\partial \alpha}=0
$$ 
$$
\frac{\partial \mathcal C}{\partial \theta}=\sum\limits_{i=1}^{n}{\left(F(x_{i:n}|\alpha,\theta)-\frac{2i-1}{2n}\right)}\frac{\partial F(x_{i:n}|\alpha,\theta)}{\partial \theta}=0
$$ 
where, $\frac{\partial F(x_{i:n}|\alpha,\theta)}{\partial \alpha}$, $\frac{\partial F(x_{i:n}|\alpha,\theta)}{\partial \theta}$ is defined in Equations ($21$), ($22$) respectively. Hence, the estimated survival and hazard rate functions are obtained as
\begin{eqnarray}\label{eq22}
\hat{S}(x)_{cme}&=& 1-\left[ 1-\left(1+{\hat{\theta}_{cme}}+{\hat{\theta}_{cme}} x+\frac{{\hat{\theta}_{cme}}^2 x^2}{2}\right)\frac{e^{-{\hat{\theta}_{cme}} x}}{1+{\hat{\theta}_{cme}}}\right]^{\hat{\alpha}_{cme}}
\end{eqnarray}
and
\begin{eqnarray}\label{eq23}
\hat{H}(x)_{cme}&=&\left\{\dfrac{\frac{\alpha\hat{\theta}_{cme}^2}{1+\hat{\theta}_{cme}}  \left[ 1-\left(1+\hat{\theta}_{cme}+\hat{\theta}_{cme} x+\frac{\hat{\theta}_{cme}^2 x^2}{2}\right)\frac{e^{-\hat{\theta}_{cme} x}}{1+\hat{\theta}_{cme}}\right]^{({\hat{\alpha}_{cme}}-1)} \left(1+\frac{\hat{\theta}_{cme} x^2}{2}\right)e^{-\hat{\theta}_{cme} x}}{1-\left[ 1-\left(1+\hat{\theta}_{cme}+\hat{\theta}_{cme} x+\frac{\hat{\theta}_{cme}^2 x^2}{2}\right)\frac{e^{-\hat{\theta}_{cme} x}}{1+\hat{\theta}_{cme}}\right]^{\hat{\alpha}_{cme}}}\right\} 
\end{eqnarray}

\subsection{Maximum product spacings estimator}
The maximum product spacing method has been introduced by Cheng and Amin ($1979,\\~1983$) as an alternative to MLE for the estimation of the unknown parameters of continuous univariate distributions. This method was also derived independently by Ranneby ($1984$) as an approximation to the Kullback-Leibler measure of information. To motivate our choice, Cheng and Amin ($1983$) proved that this method is as efficient as the MLE and consistent under more general conditions. Let us define 
\begin{equation*}
\mathcal D_{i}(\alpha,\theta)=F\left( x_{i : n}\mid \alpha,\theta \right)
-F\left( x_{ i-1 : n}\mid \alpha,\theta \right), \qquad i=1,2,\ldots ,n,
\end{equation*}
where, $F(x_{ 0 : n}\mid \alpha,\theta)=0$ and $F( x_{ n+1 : n}\mid \theta)=1-F( x_{n}\mid \theta)$ Clearly, \ $\sum_{i=1}^{n+1} \mathcal D_i (\theta) =1$. The MPSEs of the parameter $\alpha$, $\theta$ are obtained by maximizing the geometric mean of the spacings with respect to $\alpha$, $\theta$ respectively, given as
$$
GM={\left[\prod\limits_{i=1}^{n+1}\mathcal D_i(\alpha,\theta)\right]^\frac{1}{n+1}}
$$
$$
H=\log GM=\frac{1}{n+1}\sum\limits_{i=1}^{n+1} \log \mathcal D_i(\alpha,\theta)
$$
The estimates of $\alpha$ and $\theta$ are obtained by solving the following non-linear equations\\
$$\frac{\partial H}{\partial \alpha}=\frac{1}{n+1} \sum\limits_{i=1}^{n+1} \frac{1}{D_i(\alpha,\theta)} \frac{\partial D_i(\alpha,\theta)}{\partial \alpha}=0$$ 
$$\frac{\partial H}{\partial \theta}=\frac{1}{n+1} \sum\limits_{i=1}^{n+1} \frac{1}{D_i(\alpha,\theta)} \frac{\partial D_i(\alpha,\theta)}{\partial \theta}=0$$
where,
$$
\frac{\partial D_i(\alpha,\theta)}{\partial \alpha}=\frac{\partial F\left( x_{i : n}\mid \alpha,\theta \right)}{\partial \alpha}-\frac{\partial F\left( x_{ i-1 : n}\mid \alpha,\theta \right)}{\partial \alpha}
$$
$$
\frac{\partial D_i(\alpha,\theta)}{\partial \theta}=\frac{\partial F\left( x_{i : n}\mid \alpha,\theta \right)}{\partial \theta}-\frac{\partial F\left( x_{ i-1 : n}\mid \alpha,\theta \right)}{\partial \theta}
$$
 can be computed from Equations ($21$) and ($22$) respectively. Hence, the estimated survival and hazard rate functions are obtained as
\begin{eqnarray}\label{eq22}
\hat{S}(x)_{mpse}&=& 1-\left[ 1-\left(1+{\hat{\theta}_{mpse}}+{\hat{\theta}_{mpse}} x+\frac{{\hat{\theta}_{mpse}}^2 x^2}{2}\right)\frac{e^{-{\hat{\theta}_{mpse}} x}}{1+{\hat{\theta}_{mpse}}}\right]^{\hat{\alpha}_{mpse}}
\end{eqnarray}
and
\begin{eqnarray}\label{eq23}
\hat{H}(x)_{mpse}&=&\left\{\dfrac{\frac{\alpha\hat{\theta}_{mpse}^2}{1+\hat{\theta}_{mpse}}  \left[ 1-\left(1+\hat{\theta}_{mpse}+\hat{\theta}_{mpse} x+\frac{\hat{\theta}_{mpse}^2 x^2}{2}\right)\frac{e^{-\hat{\theta}_{mpse} x}}{1+\hat{\theta}_{mpse}}\right]^{({\hat{\alpha}_{mpse}}-1)} \left(1+\frac{\hat{\theta}_{mpse} x^2}{2}\right)e^{-\hat{\theta}_{mpse} x}}{1-\left[ 1-\left(1+\hat{\theta}_{mpse}+\hat{\theta}_{mpse} x+\frac{\hat{\theta}_{mpse}^2 x^2}{2}\right)\frac{e^{-\hat{\theta}_{mpse} x}}{1+\hat{\theta}_{mpse}}\right]^{\hat{\alpha}_{mpse}}}\right\} 
\end{eqnarray}


\section{Real Data Analysis}
In this Section, we have considered one medical data set which represents the survival times (in days) of $45$ patients which are observed by Gastrointestinal Tumor Study Group ($1982$) when the results of a trial comparing chemotherapy versus combined chemotherapy (CT) and radiation therapy (RT) in the treatment of locally unresectable gastric cancer. The data set is taken from Stablein and Koutron (1985). The fitting of the proposed model has been shown with the following lifetime distributions.
\begin{itemize}
	\item exponential distribution (ED), the PDF of the ED is $f(x,\theta)=\theta e^{-\theta x},~~~x,\theta >0$
	\item Lindley distribution (LD)
	\item Rayleigh Distribution (RD)
	\item xgamma distribution (XGD)
	\item Generalized exponential distribution (GED)
	\item Weibull distribution (WD)
	\item Gamma distribution (GD)
\end{itemize} 
The fitting/compatibility has been performed using different model selction tools namely, negative of log-likelihood, Akaike information criterion (AIC), corrected AIC, HQIC, BIC and KS test. The model with least AIC, CAIC, HQIC, BIC and KS is treated as best model. The obtained measures are reported in the following table which indicates that the EXGD is best choices among one parameter as well as two parameters probability distributions; hence EXGD may be chosen as an alternative model. Further, the fitted density for all the considered distribution with the histogram of the data and empirical cumulative distribution function plots are presented in Figure 3 which indicates that the proposed model provides adequate fitting to the considered data data. 
\begin{table}[htbp]
	\centering
	\caption{Values of the estimate of the parameter and model selection tools.}
	\begin{tabular}{cccccccc}\\
		\hline
		Model & MLE   &-LogL  & AIC   & CAIC  & HQIC  & BIC   & KS \\
		\hline
		ED    & 0.00139 & 340.9940 & 683.9880 & 686.7947 & 684.6616 & 685.7947 & 0.1421 \\
		LD    & 0.00278 & 342.6316 & 687.2631 & 690.0698 & 687.9366 & 689.0698 & 0.1524 \\
		RD    & 1.06E-06 & 356.6573 & 715.3147 & 718.1213 & 684.6616 & 717.1213 & 0.3321 \\
		XGD   & 0.00413 & 346.9955 & 695.9910 & 698.7976 & 696.6645 & 697.7976 & 0.2035 \\
		GED   & [1.3099, 0.00165] & 340.1633 & 684.3265 & 689.9399 & 685.6735 & 687.9399 & 0.1012 \\
		WD    & [1.1568,  0.00132] & 340.2174 & 684.4348 & 690.0481 & 685.7818 & 688.0481 & 0.4318 \\
		GD    & [1.2822,  0.00178] & 340.1868 & 684.3735 & 689.9868 & 685.7205 & 687.9868 & 0.1033 \\
		EXGD  & [0.4634,  0.00278] & 338.6061 & 681.2121 & 686.8255 & 682.5592 & 684.8255 & 0.1326 \\
	\hline
	\end{tabular}%
	\label{tab:addlabel}%
\end{table}%
\begin{center}
\begin{figure}
\includegraphics[height=4in,width=6.5in]{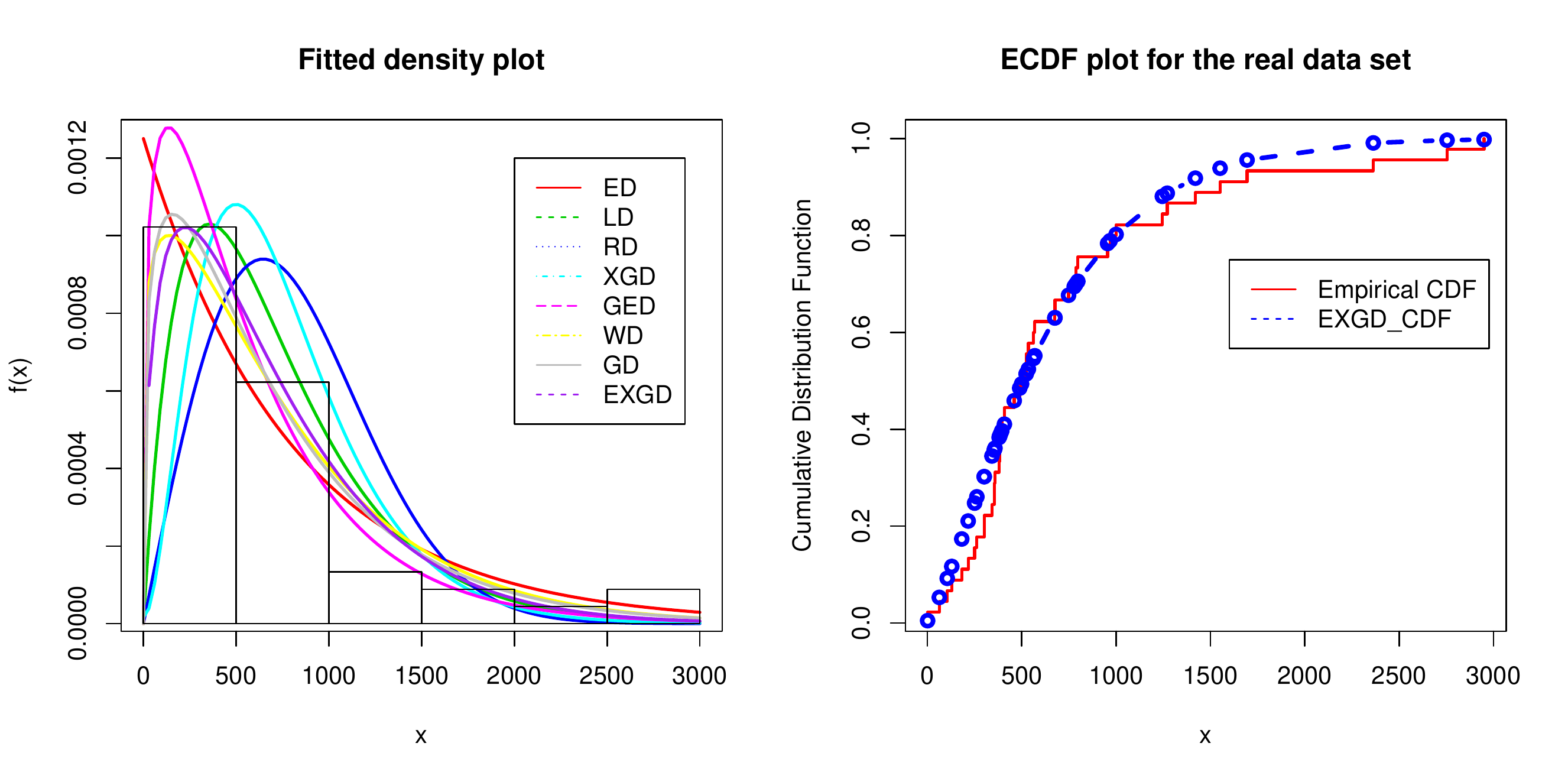}		
\caption{Fitted PDF and empirical CDF plot of EXGD.}
\end{figure}
\end{center} 
\section{Concluding remarks}
In this article, we have proposed a new positively skewed probability model, namely, EXGD by considering the generalization of XGD, introduced by Sen et al. ($2016$). Different statistical properties have been discussed. Different methods of estimation, viz., MLE, LSE, WLSE, CME and MPSE have been discussed for estimating the unknown parameters as well as reliability characteristics of the proposed model. 
Finally, a real data set has been analyzed for illustration purposes of the proposed study. Estimation of the parameters and the reliability characteristics may be further studied under different types of censoring scheme. Also the Bayesian estimation may further be considered with suitable priors and loss functions in future.\\

\end{document}